\documentclass[]{amsart}

\usepackage{amscd}
\usepackage{amsmath}
\usepackage{amsthm}
\usepackage{amsfonts}
\usepackage{amssymb}

\numberwithin{equation}{section}

\newtheorem{pro}{Proposition}

\newtheorem{lemma}[pro]{Lemma}

\newtheorem{corollary}[pro]{Corollary}
\newtheorem*{corollary*}{Corollary}
\newtheorem{thm}[pro]{Theorem}%[section]
\newtheorem*{theorem*}{Theorem}
\newtheorem*{lemma*}{Lemma}

\newtheorem*{thmA}{Theorem A}
\newtheorem*{thmB}{Theorem B}

{\theoremstyle{definition}

\newtheorem*{claim*}{Claim}

 }

\numberwithin{equation}{section}

\newcommand{\ra}{\rightarrow}

\newcommand{\Hom}{\mathop{\rm Hom}\nolimits}

\newcommand{\bb}[1]{{\mathbb #1}}    % Bourbaki 
\newcommand{\bbR}{{\bb R}}

\newcommand{\bbZ}{{\bb Z}}

\newcommand{\bbC}{{\bb C}}

\newcommand{\Alb}{\mathrm{Alb}}
\newcommand{\alb}{\mathrm{alb}}

\begin{document}
\pagestyle{myheadings}
\thispagestyle{empty}
\setcounter{page}{1}

\title[Aspherical K\"ahler Manifolds]
{Aspherical K\"ahler Manifolds with solvable fundamental group}

\author[Baues]{Oliver Baues}
\address{Mathematisches Institut II,
Universit\"at Karlsruhe,
D-76128 Karlsruhe, Germany}
\email{baues@math.uni-karlsruhe.de} 

\author[Cort\' es]{Vicente Cort\' es}
\address{Institut \'Elie Cartan, 
Universit\'e Henri Poincar\'e - Nancy I,
B.P. 239, F-54506 Van\-doeu\-vre-l\`es-Nancy Cedex, France}
\email{cortes@iecn.u-nancy.fr}

\dedicatory{F\"ur Frau K\"ahler anl\"a{\ss}lich von Erich K\"ahlers 
hundertstem Geburtstag}

\thanks{The first author gratefully acknowledges partial
support by DFG Schwerpunkt 1154, Globale Differentialgeometrie.
This work was finished at the ESI, Vienna, in October 2005. 
Both authors wish to thank ESI, and the organizers of the 
Special Research Semester on Geometry of pseudo-Riemannian manifolds
for hospitality and support.}

%\date{draft: \today}
\date{January 20, 2006}

\begin{abstract} We survey recent developments
which led to the proof of the Benson-Gordon
conjecture on K\"ahler quotients of
solvable Lie groups. In addition we 
prove that the Albanese morphism of a 
K\"ahler manifold which is a homotopy torus is a
biholomorphic map. The latter result then implies
the classification of compact aspherical K\"ahler 
manifolds with (virtually) solvable fundamental group
up to biholomorphic equivalence. They are all
biholomorphic to complex manifolds which
are obtained as a quotient of $\bbC^n$ by
a discrete group of complex isometries.
%
%that every aspherical K\"ahler 
%manifolds with solvable by finite fundamental group
%is biholomorphic to complex a 

%
%and gives their diffeomorphism classification
%This proves the K\"ahler rigidity
%of smooth aspherical manifolds
%with solvable fundamental group
%and gives their diffeomorphism classification.  
%This result implies the classification of
%aspherical manifolds with K\"ahler structure
% 
\end{abstract}

\maketitle

\section*{}
The purpose of this article is to give an account on 
recent developments concerning the classification of
compact aspherical K\"ahler manifolds whose
fundamental groups contain a solvable subgroup 
of finite index. This problem stirred some interest
among differential geometers for quite some time.
It traces back to the classification of K\"ahler 
Lie groups and became popular through a  
conjecture of Benson and Gordon on K\"ahler quotients 
of solvable Lie groups. 
Only recently satisfying answers to the 
original Benson-Gordon
conjecture became available when 
methods from complex geometry were introduced
into the subject.  
%  when complex geometers took interest in the problem. 
The story as it presents itself now may be seen as an interesting
show case where research questions from 
differential geometry, complex geometry and topology meet. 

%
%Hence, one particular aim  of our paper is to make the various 
%aspects of the subject accessible to
%mathematicians working in either field. (in a broad spectrum....) 

\subsubsection*{Outline of the paper}

The first section of the paper
has a survey character. We introduce
the Benson-Gordon 
conjecture and describe its relation to the 
classification of K\"ahler Lie groups, K\"ahler groups
and K\"ahler manifolds.  
%We explain this 
%in the general context 
%of aspherical manifolds with solvable
%fundamental group. 
In the course we also 
sketch the proof of a theorem of 
Arapura and Nori on polycyclic K\"ahler
groups. This leads us to the following 
rigidity result for K\"ahler infra-solvmanifolds, 
which, in particular, contains the solution 
to the conjecture of Benson and Gordon.
  
\begin{thmA} Let $M$ be an infra-solvmanifold which admits
a K\"ahler metric. Then $M$ is diffeomorphic to a flat 
Riemannian manifold.  
\end{thmA}  

In section 2 we develop generalisations of the 
previously described results. These concern 
the classification of K\"ahler structures on 
general smooth aspherical manifolds with polycyclic 
fundamental group. One new key point is to establish
that a K\"ahler manifold which is homotopy equivalent
to a torus is also diffeomorphic to a torus. 
We prove the stronger result that a K\"ahler manifold which is 
a homotopy torus is biholomorphic to a complex torus. 
This result is based on a detailed study of the
Albanese mapping of such K\"ahler manifolds, which 
is  carried out in section 3.  In section 2.2 we 
introduce our main result:

\begin{thmB} Let $X$ be a compact 
complex manifold which is aspherical 
with virtually solvable fundamental group and assume 
$X$ supports a K\"ahler metric.
Then $X$ is biholomorphic to a quotient of $\bbC^n$ by
a discrete group of complex isometries.
\end{thmB}

The theorem shows that the existence of a K\"ahler metric 
on a complex aspherical manifold with virtually solvable fundamental group 
forces not only the the fundamental group, but also the smooth structure 
and the complex structure to be as simple as possible. \\

\noindent 
In order to provide the reader with an overall view,  
we list the main arguments in the proof of Theorem B, which 
are essentially the following:\\ 
\begin{enumerate}
\item[(i)] A solvable group  
which is the fundamental group of  
a compact {\em aspherical} manifold is torsion-free and 
polycyclic (by results of Bieri mentioned in 1.5).   
\item[(ii)] A polycyclic K\"ahler group
is virtually nilpotent (by a theorem of Arapura and Nori stated in 1.5 and 
explained  in 1.6).  
\item[(iii)] A nilpotent \emph{aspherical} K\"ahler group is abelian 
(by work of Benson and Gordon, Hasegawa, see 1.3). 
\item[(iv)] A K\"ahler manifold which is homotopy equivalent
to a complex torus is biholomorphic to a complex torus (see Theorem 
\ref{Thm2} and Theorem \ref{Thm5} in 2.2, section 3 for the proofs). 
\end{enumerate}  

%We hope that our article reveals some of 
%the rather subtle problems which are involved in
%this subject. 

%
\subsubsection*{Acknowledgement} 
We thank J\"org Winkelmann for discussions which were 
helpful for the proof of Theorem 5. 
We are grateful to Fr\'ederic Campana for his detailed 
comments on the first draft of our paper. 
We also wish to thank Frank Herrlich for further discussions,
and Vincent Koziarz for providing us reference \cite{Lelong}.

\section{From K\"ahler Lie groups to solvable
K\"ahler groups} % and aspherical manifolds}
Recall the differential geometer's definition of a K\"ahler manifold. 
A \emph{K\"ahler manifold} is a Riemannian 
manifold endowed with a parallel skew-symmetric complex structure,
see \cite{KN_2}. The closed non-degenerate two-form which is associated
with such data is called
the K\"ahler form. Notice that, by the Newlander-Nirenberg theorem, 
a K\"ahler manifold has an underlying structure of a 
complex manifold. Our paper is concerned with the existence
problem for K\"ahler structures on solvmanifolds and related
spaces. The story starts, maybe, with the study  
of homogeneous K\"ahler manifolds and, in particular, 
K\"ahler Lie groups. 

\subsection{K\"ahler Lie groups}
A Lie group is called K\"ahler if it is equipped with a left-invariant
K\"ahler metric and a left-invariant K\"ahler form.  
The basic examples are \emph{flat K\"ahler Lie groups} and solvable 
K\"ahler Lie groups acting simply transitively on a bounded 
domain in $\bbC^n$ by biholomorphisms \cite{GPSV}. 
There is a fundamental conjecture concerning the classification 
of homogeneous  K\"ahler manifolds up to biholomorphism \cite{GPSV}, which 
was finally settled in \cite{DN}. As a reference for the rich 
structure of solvable K\"ahler Lie groups and their classification  
up to (holomorphic) isometry one may consult \cite{PS}.  

Typically, a solvable Lie group endowed with a left-invariant
complex structure can admit many non-equivalent irreducible 
K\"ahler Lie group structures 
%Let us point out that, in the latter 
%case,  the invariant K\"ahler metric is not necessarily a multiple of 
%the Bergmann metric  and that a K\"ahler Lie group is not necessarily 
%a Riemannian product of K\"ahler-Einstein Lie groups, 
(contrary to what is claimed in \cite{Li}). 
However, \emph{a unimodular 
K\"ahler Lie group $G$ is always flat}, as proved by Hano \cite{Hano}.  
In particular, the Lie group $G$ is solvable of a very restricted type.
This follows from the classification of flat Riemannian 
Lie groups which is given in \cite{Milnor_2}. 
As a particular consequence, a K\"ahler Lie group $G$ which admits a 
lattice, that is, a discrete cocompact subgroup 
$\Gamma \leq G$, is flat. Hence, any compact quotient space $G\big/\Gamma$,
where $G$ is a 
K\"ahler Lie group, is a compact flat Riemannian manifold, and,  
by Bieberbach's theorem 
\cite{Wolf}, it  is finitely covered by a flat torus.
 
Let $G$ be a solvable Lie group which admits a lattice $\Gamma$. 
Independently of the existence of a K\"ahler Lie group structure on $G$, 
the question arises whether the compact smooth manifold $G\big/\Gamma$ 
carries a K\"ahler structure. Benson and Gordon initiated the
study of this question in \cite{Be-Go_1, Be-Go_2}. 
Their conjecture,  which we explain in the next paragraph, suggested
an answer to this question for quotient spaces of completely 
solvable Lie groups. 

\subsection{A conjecture of Benson and Gordon}
A solvable Lie group is called \emph{completely solvable} if its 
adjoint representation  has only real eigenvalues. 
In particular,
nilpotent Lie groups are completely solvable.
Notice that a completely 
solvable unimodular K\"ahler Lie group
is abelian in virtue of Hano's result and 
the classification of flat Lie groups as contained 
in \cite{Milnor_2}. 

Lattices in 
completely solvable groups satisfy 
strong rigidity properties, analogous to the Malcev
rigidity  \cite{Malcev} of lattices in nilpotent Lie groups.
This was first proved by Saito, 
see \cite{Saito}. In particular, the embeddings of $\Gamma$ in $G$ are 
conjugate by an automorphism of $G$, and the completely 
solvable Lie group $G$ is determined by $\Gamma$
up to isomorphism. These results show a particular strong link
between the properties of the fundamental  group of the 
space $G\big/\Gamma$ and the group $G$.
This may serve as a motivation for the
following 

\subsubsection*{Benson-Gordon conjecture \cite{Be-Go_1, Be-Go_2}:}
Let $M= G\big/\Gamma$ be a compact quotient
space of a completely solvable
Lie group $G$ by a discrete subgroup $\Gamma$. 
If $M$ admits a  K\"ahler metric then $M$ is \emph{diffeomorphic} 
to a standard $2n$-torus $T^{2n} = S^1 \times \cdots \times S^1$.\\

%Moreover, as proved
%recently by Wilking \cite{Wilking}, any isometric affine action on  
%a solvable Lie group is isometric to an affine action
%on a solvable Lie group.  

\subsection{K\"ahler nilmanifolds}

A nilmanifold $N\big/\Gamma$ is a quotient space of a simply 
connected nilpotent
Lie group $N$ by a discrete cocompact subgroup $\Gamma \leq N$.
Nilmanifolds satisfy the
Benson-Gordon conjecture. This means that
\emph{any 
compact nilmanifold which is K\"ahler is
diffeomorphic to a torus.} 

The reasons for this fact
stem from the Hodge theory of K\"ahler manifolds.
The de Rham cohomology ring of a
compact K\"ahler manifold satisfies strong extra conditions,
which arise from the existence of the K\"ahler class and from 
Hodge theory. Elementary restrictions are that the second Betti-number
is positive and odd Betti numbers are even. 
See \cite{Griffiths-Harris, Voisin, Wells} 
for basic reference on this subject, and \cite{DGMS}
for the formality of the minimal model of a K\"ahler manifold. 

Now by a result of Nomizu \cite{Nomizu}, the de Rham 
cohomology of a nilmanifold $N\big/\Gamma$ is
computed by the finite-dimensional Koszul complex
for the Lie algebra cohomology of $N$. It was observed
in \cite{Be-Go_1} that if the Lie algebra cohomology 
of $N$ satisfies the \emph{strong Lefschetz property}, 
which holds for the cohomology ring of a K\"ahler manifold, 
then $N$ is abelian. Also the Koszul complex of $N$ is
a minimal model (in the sense of Sullivan, see \cite{DGMS})
of the de Rham complex for $N\big/\Gamma$.
The minimal model of a K\"ahler manifold
is formal \cite{DGMS}. As observed in \cite{Hase_1}, if the 
Lie algebra cohomology of $N$ is
\emph{formal}, $N$ must be abelian.
Another approach, announced in  \cite{CFG},
uses the fact that Massey products in 
K\"ahler manifolds vanish \cite{DGMS}. 

However, these techniques did not suffice to prove the
full conjecture of Benson and Gordon. (Note that the proof given
in \cite{Tralle-Kedra} contains a mistaken 
argument.) Moreover, there do
exist compact quotient spaces of completely solvable Lie groups with
non virtually nilpotent fundamental group, and
which have the cohomology ring 
of a K\"ahler manifold, see \cite{Be-Go_2,Sa-Ya}.

\subsection{Aspherical  K\"ahler groups}

In general, it is a delicate problem to understand which finitely presented 
groups can occur as the fundamental group of a compact K\"ahler manifold. 
See \cite{ABC} for some results and techniques. Traditionally, 
the fundamental group of a compact 
K\"ahler manifold is called a \emph{K\"ahler group}.
A manifold is called 
aspherical if its universal covering
space is contractible. We call a group an \emph{aspherical K\"ahler group} 
if it is the fundamental group of a compact aspherical
K\"ahler manifold. 

The homotopy properties of an aspherical 
manifold, its \emph{homotopy type}, and, in particular,
the cohomology
ring and minimal model depend only on its
fundamental group. Recall that every finitely generated 
torsion-free nilpotent group occurs as the fundamental 
group of a compact nilmanifold \cite{Malcev}. Thus,
the results of 
Benson and Gordon \cite{Be-Go_1}, Hasegawa 
\cite{Hase_1} show that \emph{a
finitely generated torsion-free nilpotent 
group is not an aspherical K\"ahler group, 
unless $\Gamma$ is abelian.} 

Note that this fact, rephrased in topological terms,
also shows that \emph{a compact aspherical K\"ahler manifold with
nilpotent fundamental group is 
homotopy equivalent to a torus $T^{2n}$, 
for some $n$}. 

Let us further remark that, although there are no 
non-abelian nilpotent \emph{aspherical}  K\"ahler 
groups, there are (as shown by Campana \cite{Campana}) 
examples of torsion-free two-step nilpotent 
(in particular not virtually-abelian) K\"ahler groups.  
Despite these pioneering examples,  nilpotent  K\"ahler groups 
remain a complete mystery.  

\subsection{Polycyclic K\"ahler groups}
A \emph{polycyclic group} is a group which
admits a finite normal series with
cyclic quotients. See \cite{Segal} for
a comprehensive account on the theory 
of these groups. Note that, in particular,
finitely generated nilpotent groups
are polycyclic. But the class of polycyclic groups 
is much larger. By a result of Mostow \cite{Mostow_3},
the fundamental group of a quotient space of a solvable 
Lie group is a polycyclic group. More generally, a result
of Bieri \cite{Bieri, Bieri_2} characterizes torsion-free
polycyclic groups as solvable Poincar\'e duality
groups.  This implies that \emph{a solvable group 
which is the fundamental group of
a compact aspherical manifold is 
a torsion-free polycyclic group}.

Let $\mathcal P$ be a property of groups. Recall that a group is
said to be virtually  ${\mathcal P}$ if it has a finite
index subgroup which satisfies ${\mathcal P}$.
Concerning polycyclic K\"ahler groups there is the following remarkable result
which is proved in \cite{Arapura-Nori}: 

\begin{theorem*}[Arapura-Nori] 
Let $\Gamma$ be
a (virtually) polycyclic K\"ahler group
then $\Gamma$ is virtually nilpotent.
\end{theorem*}

We shall explain the proof 
of this key result further below. The result of
Arapura and Nori as stated in \cite{Arapura-Nori} 
is more general and concerns also quasi-K\"ahler groups. 
It was shown by Campana \cite{Campana_2} that the assumption 
``polycyclic'' in the above theorem can be generalised 
to ``linear solvable''. However, it is not known whether there exist 
(non-linear) solvable K\"ahler groups which are not virtually nilpotent.

Combining the theorem with the 
results mentioned in 1.4 and 1.5 yields:

\begin{corollary*} Every aspherical 
solvable K\"ahler group $\Gamma$ contains
a finitely generated abelian subgroup 
of finite index. In particular, $\Gamma$
is a Bieberbach group. 
\end{corollary*}

Recall that, by definition,  a Bieberbach group is
a finitely generated torsion-free group
which contains an abelian subgroup 
of finite index. A version of the
corollary is also mentioned in \cite{Mic}. 

The Benson-Gordon conjecture is actually a 
straightforward consequence of the corollary. 
We shall explain this in 1.6. 
The relevance of the 
Theorem of Arapura and
Nori for the Benson-Gordon
conjecture was first recognized by 
Arapura in \cite{Arapura_2}.  His main
result is concerned with the biholomorphic
properties of K\"ahler solvmanifolds and a
conjecture of Hasegawa \cite{Hase_2}. 
(Note, however, that the proof
contains a gap.  See MR 2038772 (2005a:53123).) 

Before continuing with the Benson-Gordon conjecture
we shall briefly describe the
techniques which are behind the proof of
the above theorem of Arapura and Nori.
The  proof is based on ideas 
and results of Beauville \cite{Beauville_1, Beauville_2} 
on the character varieties of K\"ahler groups.

\subsection{The character variety of a polycyclic K\"ahler group}
Let $\Gamma$ be a finitely generated group, and let be $k$ the
number of elements of a fixed generating set. 
Then, as is easy to see, the set of all homomorphisms
from $\Gamma$ to $\bbC^*$ forms an affine algebraic
variety, the \emph{character variety} $X_\Gamma \subset (\bbC^*)^k$ 
of $\Gamma$. Each $\rho \in X_\Gamma$ defines a one-dimensional $\Gamma$-module  $\bbC_\rho$. 
Let $X^1_\Gamma$ denote the subset of 
characters of $\Gamma$ which satisfy 
$\dim H^1(\Gamma, \bbC_\rho) >0$. 
%Then $X^1_\Gamma$
%is an algebraic subset of $H^1(\Gamma, \bbC_\rho)$. 

Now let $\Gamma$ be a K\"ahler group. Note that,  for each 
$\rho \in X_\Gamma$, the $\Gamma$-module  $\bbC_\rho$ 
defines a flat holomorphic line bundle $L_\rho$ over
any K\"ahler manifold $X$ with $\pi_1(X) = \Gamma$. 
Moreover, $H^1(\Gamma, \bbC_{\rho}) = H^1(X,L_{\rho})$.
Using this correspondence, geometric arguments and (non-abelian)
Hodge theory on $X$ provide
information about the set  $X^1_\Gamma$.
According to Beauville's main result \cite[Corollaire 3.6]{Beauville_2},
$X^1_\Gamma$ is a finite union of translates of connected subgroups 
of $X_{\Gamma}$ (which arise from homomorphisms of $\Gamma$
to surface groups) and a finite number of unitary characters. 
In particular,  \emph{isolated points  of $X^1_{\Gamma}$ are unitary 
characters}. Moreover, it was conjectured in \cite{Beauville_2} that the 
isolated points of $X^1_\Gamma$ are characters of finite order. 
This conjecture was proven for the fundamental groups of smooth projective 
algebraic varieties by Simpson \cite{Simpson} and for general K\"ahler groups 
by Campana \cite{Campana_2}. 
For polycyclic K\"ahler groups the conjecture 
follows from a direct argument 
which is based on the next lemma. 

% This follows from the following lemma which is contained in 
%\cite{Beauville_2}: 

\begin{lemma*} Let $\Gamma$ be a finitely generated group
with $D \Gamma\big/ D^2 \Gamma$ finitely generated.
Then $X^1_\Gamma$ is a finite set. If, in addition, $\Gamma$ 
is a K\"ahler group, then the characters of $X^1_\Gamma$
are of finite order. 
\end{lemma*}

Here,  $D^i \Gamma$ denotes the derived series of 
$\Gamma$. We briefly recall the proof of this Lemma, 
which is given in \cite{Beauville_2}. 

\begin{proof} 
The formula 
$H^1(\Gamma, \bbC_\rho)  = 
\Hom_{\Gamma^{\rm ab}}(D\Gamma/D^2 \Gamma, \bbC_\rho)$,
where $\Gamma^{\rm ab}= \Gamma\big/D\Gamma $ 
denotes the abelianisation,  shows that the set $X^1_\Gamma$ is finite. 
Moreover, the values of any $\rho \in X^1_\Gamma$  
on $\Gamma$ are algebraic integers, since, as the formula
shows, $\rho$ occurs as a weight 
of the adjoint representation of $\Gamma$ on 
$D\Gamma/D^2 \Gamma \otimes \bbC$.
If $\Gamma$ is K\"ahler
then, by the structure result on $X^1_{\Gamma}$, $\rho$ 
is unitary. Therefore, these algebraic integers are of absolute value one,
and also all of their algebraic conjugates are.
Hence, by Kronecker's Lemma, these algebraic
integers are roots of unity. It follows that
$\rho$ is of finite order. 
\end{proof}

A polycyclic group satisfies the assumption
of the lemma. In fact, every subgroup of $\Gamma$ is finitely 
generated, see \cite{Segal}. Hence, \emph{for a polycyclic K\"ahler group, 
the elements of $X^1_\Gamma$ are of finite order}.
As observed in \cite{Arapura-Nori}, we furthermore have:

\begin{lemma*}  Let $\Gamma$ be a polycyclic group with  $X^1_\Gamma$ 
consisting of characters of finite order, then $\Gamma$ is virtually nilpotent.
\end{lemma*}

A direct proof for this remark (following the lines of
\cite{Arapura-Nori}) can be given, by considering the embedding of 
$\Gamma$ into its algebraic hull $H_{\Gamma}$. The hull $H_{\Gamma}$
is a linear algebraic group which contains
$\Gamma$ as a Zariski-dense subgroup, and furthermore 
it has the property that
the centraliser of the unipotent radical $U$ of
$H_{\Gamma}$ is contained in $U$.  (See \cite{Borel}, for
the general concept of a linear algebraic group, 
and \cite{Mostow_2}, \cite{Arapura-Nori} and also \cite[Appendix A]{Baues} 
for details on the construction of $H_{\Gamma}$ and 
further applications.) 

\begin{proof}
By replacing $\Gamma$ with a finite index subgroup  
we can assume that $H_{\Gamma}$ is connected.
By standard theory of linear algebraic groups \cite{Borel}, the  
group $H_{\Gamma}$ has an algebraic splitting
$H_{\Gamma} = T U$, where $T$ is an algebraic torus,
that is, $T$ is a connected abelian group of semi-simple elements.
Note that $D \Gamma \subset U$, and  $N= D H_{\Gamma}$ is the 
Zariski-closure of $D \Gamma$.
Since the centraliser of $U$ is contained in $U$,
the abelian reductive group $T$ is faithfully represented
on $N$, and therefore (as holds for reductive actions
on unipotent groups)  also on $N/D N = D\Gamma/D^2 \Gamma \otimes \bbC$. 
The weights of $T$ on $N/D N$ restrict to elements of $X^1_{\Gamma}$.
Since $\Gamma$ is
Zariski-dense, also the weights of $T$ are of finite order. 
Since $T$ is connected, they are trivial. Therefore, $T$
centralises $U$.
This implies that $H_{\Gamma}=U$ is a nilpotent group,  and
so is $\Gamma$.  
\end{proof}

\subsection{Proof of the  Benson-Gordon conjecture} 

The above corollary to the theorem of Arapura and Nori
easily implies
the original Benson-Gordon 
conjecture.  In fact, let $G$ be simply connected 
and completely solvable, and let $\Gamma \leq G$ a lattice.
By the
corollary, the fundamental group
$\Gamma$ of a compact K\"ahler
manifold $G\big/\Gamma$ is virtually abelian. 
The Saito rigidity for lattices in completely
solvable Lie groups then implies
that $G$ and, therefore, also $\Gamma$
are abelian, and moreover, $G\big/\Gamma$
is diffeomorphic to the torus 
$T^{2n}$.

\subsection{K\"ahler infra-\-solv\-manifolds}
In the application of the Corollary in 1.5, 
it is not necessary to restrict oneself to the
case of quotient spaces of completely solvable
Lie groups. The quotient space $G/H$ of a solvable 
Lie group $G$ by a closed subgroup $H$ 
is called a \emph{solvmanifold}.
More generally, we can consider \emph{infra-solvmanifolds}. 
Infra-solvmanifolds are constructed using isometric
affine actions on solvable Lie groups, see
\cite{Baues, Wilking} for precise definition  
and further reference.  
This class of manifolds contains, in particular,  
solvmanifolds and flat Riemannian manifolds. 
As is well known, every Bieberbach group 
is the fundamental group of a compact 
flat Riemannian manifold, see \cite{Wolf}. 
(Also the converse is true.) 
Moreover, see \cite{Baues}, \emph{every  torsion-free
virtually polycyclic group occurs as the fundamental
group of a compact infra-solvmanifold}.

Flat Riemannian manifolds with isomorphic fundamental groups 
are diffeomorphic, by Bieberbach's theorem \cite{Bieberbach}.
The same is true for compact solvmanifolds, by Mostow's 
famous result \cite{Mostow_1}. Only recently (see
\cite{Baues}), it was 
established that infra-solvmanifolds are 
\emph{smoothly rigid}, that is, these manifolds
are determined by their fundamental group
up to diffeomorphism. Therefore, 
infra-solvmanifolds constitute natural \emph{smooth model spaces}
for compact aspherical manifolds with 
a torsion-free virtually polycyclic 
fundamental group. 
The smooth rigidity of
infra-solvmanifolds and
the Corollary in 1.5 
imply the following --- generalised --- 
version of the Benson-Gordon 
conjecture:

\begin{corollary*} Let $M$ be an (infra-)
solvmanifold which admits
a K\"ahler metric. Then $M$ is
diffeomorphic to a flat 
Riemannian manifold.  
\end{corollary*}

But, of course, the story does not end here.   

\section{Classification of aspherical K\"ahler manifolds}

We start this section with a short digression
on the existence of exotic smooth structures on
aspherical manifolds. Then we discuss the classification 
of aspherical K\"ahler manifolds with a virtually 
polycyclic fundamental group up to diffeomorphism
and up to biholomorphic equivalence. By the previously
discussed results on polycyclic K\"ahler groups, 
these classifications are reduced 
to the classification of  K\"ahler structures on
homotopy tori and their quotient spaces. The latter
classification is obtained by studying the Albanese mapping of 
such a K\"ahler manifold.

\subsection{Diffeomorphism classification}

\subsubsection*{Fake tori}
Let $T^n = S^1 \times \cdots \times S^1$ denote the
$n$-torus with its  natural differentiable structure. 
A compact topological manifold is 
called a \emph{homotopy torus} if it is homotopy 
equivalent to $T^n$. For the classification of homotopy tori up
to piece-wise linear equivalence, see \cite{Hsiang_Shane}. 
Every homotopy torus is smoothable, but not uniquely.   
By the work of Browder \cite{Browder} and 
Wall \cite[\S 15 A, B]{Wall}, Hsiang and Shaneson \cite{Hsiang_Shane}
there do exist $n$-dimensional smooth homotopy-tori, $n \geq 5$,
which are  \emph{not} diffeomorphic to $T^n$. 
One way to obtain such \emph{fake tori}
is to glue with exotic spheres \cite{Browder}. 
By the classification of PL-structures on $T^n$, 
see \cite{Hsiang_Shane},
there exist even more fake tori. Furthermore, 
all fake tori admit a finite
smooth covering by a standard torus 
\cite[\S 15 A, B]{Wall}, and are 
homeomorphic to $T^n$ \cite{Hsiang_Wall}.

Recall from 1.6 that, by  Mostow's theorem
\cite{Mostow_1}, the natural smooth structure
on a solvmanifold $M=G/H$, which is inherited from $G$ and
$H$, is independent of the choice of $G$ and $H$,
and it is determined by the fundamental group of $G/H$
only. As already mentioned, the analogous fact also 
holds for the larger class
of  infra-solvmanifolds. The results
cited above show that, surprisingly, there do exist smooth
manifolds, even finite smooth quotients of a standard torus $T^n$,
which are homeomorphic but not diffeomorphic to an (infra-) 
solvmanifold. 

\subsubsection*{Aspherical K\"ahler manifolds with solvable fundamental group}
In the light of these facts, it seems natural to ask for the
diffeomorphism classification
of \emph{all} compact aspherical
K\"ahler manifolds with virtually solvable fundamental group, 
instead of restricting our interest to quotient 
manifolds of solvable Lie groups only. The answer to
this question  is contained in the following generalisation
of the result in 1.8:

\begin{thm} Let $M$ be a
compact aspherical K\"ahler 
manifold with virtually 
solvable fundamental group. 
Then $M$ is diffeomorphic to a 
flat Riemannian 
manifold. 
\end{thm}

This result is implied by the
classification of K\"ahler structures
up to biholomorphism 
which is our next topic. 

\subsection{Classification up to biholomorphism}
Recall that a K\"ahler
manifold is equipped with a unique underlying 
complex manifold structure. Obvious obstructions for a complex structure
to be K\"ahler arise from Hodge-theory, see
\cite{Voisin}. Closely related, 
to the problem of classifying 
K\"ahler manifolds up to diffeomorphism,
is the finer classification problem for the 
corresponding underlying complex manifolds.
 
%There are many examples of complex manifold 
%structures on the standard torus $T^{2n}$. 
%Not every complex 
%structure on $T^{2n}$ is the underlying complex
%structure of a K\"ahler metric, as the classical examples 
%of ???? show. 
 A \emph{complex torus} is 
a complex manifold $\bbC^n\big/\, \Gamma$, where
$\Gamma \cong \bbZ^{2n}$ is a lattice in $\bbC^n$.
Evidently, complex tori support a flat K\"ahler metric
which is inherited from the flat K\"ahler metric on $\bbC^n$.

Our principal result in this section asserts that complex tori 
are the only examples of complex manifold structures 
on homotopy tori which can be K\"ahler: 

\begin{thm} \label{Thm2} 
Let $X$ be a compact complex manifold which
is K\"ahler and homotopy equivalent to a torus $T^{2n}$.  
Then $X$ is biholomorphic to a complex torus.  
\end{thm}

The theorem is a direct consequence of Theorem 5 below.
As F.\ Catanese kindly informed us, Theorem 2 may
also be deduced as a special case of  \cite[Proposition 4.8]{Catanese}. \\

Together with the results of Bieri \cite{Bieri},  Benson-Gordon \cite{Be-Go_1}, Hasegawa 
\cite{Hase_1}, Arapura-Nori \cite{Arapura-Nori}, 
Theorem 2 implies: 

\begin{corollary} Let $X$ be a compact 
complex manifold which is aspherical 
with virtually solvable fundamental group and assume 
$X$ supports a K\"ahler metric.
Then $X$ is biholomorphic to a finite quotient of a complex torus.  
\end{corollary}

We obtain the classification of K\"ahler aspherical
manifolds with solvable fundamental group up to biholomorphism, 
as formulated in Theorem B in the introduction:
%which generalises Arapuras classification of K\"ahler solvmanifolds 
%\cite{Arapura_2}, 
%as follows:

\begin{corollary} \label{mainCor} Let $X$ be a compact 
complex manifold which is aspherical 
with virtually solvable fundamental group and assume 
$X$ supports a K\"ahler metric.
Then $X$ is biholomorphic to a quotient of $\bbC^n$ by
a discrete group of complex isometries.
\end{corollary}

\begin{proof} By Corollary 3, $X$ is finitely covered 
by a complex torus, and the deck transformation group is
a finite group of holomorphic maps of a complex torus. 
Every biholomorphic map of a complex torus is induced
by a complex affine transformation of the vector space 
$\bbC^n$, see \cite{Griffiths-Harris}.
Thus it follows that $X$ is a quotient space
of $\bbC^n$ by a discrete group $\Gamma$ of complex affine transformations.
Since the translations in $\Gamma$ form a subgroup of finite index, 
$\Gamma$ preserves an invariant scalar product on $\bbC^n$.
\end{proof}
Since the compact complex manifolds occuring in
Corollary 4 are quotient spaces of $\bbC^n$ by affine isometries,
these manifolds also carry a flat Riemannian metric. 
Therefore, Corollary 4 implies Theorem 1.

For a discussion of the classification of discrete groups of
complex isometries on $\bbC^n$, as appearing in Corollary 4,
in further detail, see also \cite{Mic}.   

\subsubsection*{The Albanese torus of a K\"ahler manifold}
We shall invoke now a natural construction
from complex geometry which is  the main tool 
in our proof of
Theorem 2. Let us denote with $X$ a 
complex manifold, and let us put $\dim X$ for its
complex dimension. The \emph{Albanese torus} $\Alb(X)$ of $X$ 
is a complex torus  which
is naturally associated to any compact K\"ahler manifold $X$. 
Its complex dimension is $\dim \Alb(X) = n$, where $2n =\dim H^1(X,\bbC)$.
See Section 2.1, and \cite{Blanchard, Campana, Voisin}, 
or \cite[Chapter 1]{ABC} for further notes on this topic.

Roughly, the Albanese torus is constructed as follows: 
Let $H^0(X, \Omega^1_{X})$
denote the complex vector space of holomorphic one-forms on $X$. 
Since $X$ is K\"ahler, Hodge theory implies that
$\Alb(X) = H^0(X, \Omega^1_{X})^*\big/H_{1}(X,\bbZ)$ is a
complex torus. There is a natural holomorphic map, the
\emph{Albanese morphism}, $\alb_{X}: X \ra \Alb(X)$, where
$\alb_{X}(x) = \int_{x_{0}}^x$ associates to $x$ the 
integral along any path from (a fixed) $x_{0}$ to $x$, acting
on holomorphic one-forms. 

Theorem 2 is a consequence of the following rigidity result
for K\"ahler manifolds. Roughly speaking, the result asserts that 
a K\"ahler manifold which has the integral cohomology ring of a torus $T^{2n}$  
must be biholomorphic to a complex torus: 

\begin{thm}\label{Thm5}  
Let $X$ be a compact  complex manifold, $n = \dim X$,
which is K\"ahler and satisfies $2 \dim X = \dim H^1(X,\bbC)$.  If furthermore 
$\dim H^2(X,\bbC) \leq {2n \choose 2}$ and 
$H^{2n}(X,\bbZ)$ 
is generated by integral classes of degree 1, 
then the Albanese morphism $\alb_{X}: X \ra \Alb(X)$ is 
a biholomorphic map.
In particular, $X$ is biholomorphic to a complex torus.  
\end{thm}

A version of this theorem under the (stronger) assumption that 
$X$ is diffeomorphic to a $2n$-torus is stated in \cite{Arapura_2}
without proof.

\subsubsection*{Projective examples, abelian varieties}
Variants of our results also hold in the 
algebraic setting.  A smooth projective variety
is a compact submanifold of complex projective
space which is defined by algebraic equations.   
Since projective space is K\"ahler, smooth
projective varieties are naturally K\"ahler manifolds as well.
An \emph{abelian variety} is an irreducible complex
projective variety with an abelian group
structure. One can show (cf.\ \cite[\S 12.1.3]{Voisin}), 
if $X$ is projective, then $\Alb(X)$ is projective and
an abelian variety. Furthermore, the map $\alb_{X}$ is
a morphism of algebraic varieties. 

\begin{corollary} Let $X$ be a smooth projective variety 
which is homotopy-equivalent to a torus $T^{2n}$ then 
$X$ is an abelian variety.   
\end{corollary}
\begin{proof} By the above, $\alb_{X}$ is a biholomorphic 
morphism $X \ra \Alb(X)$ of projective varieties. By the 
GAGA theorem,  $\alb_{X}$ is actually an algebraic isomorphism,
that is, the holomorphic inverse of $\alb_{X}$
is a morphism of algebraic varieties,
cf.\ \cite[Appendix B, Exerc.\ 6.6]{Hartshorne}.
\end{proof}

Remark that, in complex dimensions $n=1$ and $n=2$,
all complex structures on $T^{2n}$ are K\"ahler. In 
dimension $n=1$, this fact is elementary,  for $n=2$, it
may be deduced from the classification of surfaces, as
described in \cite[VI 1.]{BPV}. In particular, by Corollary 4, 
\emph{any}  complex surface  $X$ homotopy equivalent to $T^4$
is biholomorphic to a complex torus.  
However, the analogue of this 
result fails in higher dimensions. In fact, 
Sommese  \cite[Remark III, E)]{Sommese} observed the
existence of exotic non-K\"ahler complex structures 
on the 6-dimensional real torus. These examples are
also discussed in Catanese's paper \cite{Catanese}, where the
classification problem for complex structures on 
$K(\pi,1)$ compact manifolds is formulated in a spirit
similar to the point of view of our paper. \\

We conclude our discussion with yet another research question. 
Though,  as is clear now,  nil- and solvmanifolds are \emph{not} K\"ahler, 
unless they are tori, many of them admit a complex manifold structure,
see, for example, \cite{Borcea,GMPP, Con_Fin,Tralle}.
It seems tempting to assert that the existence 
of a complex manifold structure on a smooth manifold which is
\emph{homotopy equivalent to a solvmanifold} forces
the differentiable structure to be as nice as possible, that
is, such a manifold should be diffeomorphic to a solvmanifold. Hence, 
\emph{is a homotopy torus with complex manifold structure
diffeomorphic to a standard torus?}  More generally, we ask: 
\emph{Let $X$ be a compact complex manifold, aspherical with 
virtually solvable fundamental group. 
Is $X$ diffeomorphic to an infra-solvmanifold?}\\

The rest of our paper is
devoted to give a  proof of Theorem 5.

\section{The Albanese morphism of a K\"ahler manifold}

% In this section we give the proof of Theorem 5. 
We start by defining
the Albanese variety $\Alb(X)$ of a K\"ahler manifold $X$ and
deduce some elementary properties of the Albanese
morphism $\alb_{X}$.  For further reference on $\Alb(X)$
one may consult  \cite{Blanchard, Campana}, \cite[Chapter 1]{ABC}
and \cite[\S 12.1.3]{Voisin}.
 
\subsection{The Albanese morphism}
Let $X$ be compact K\"ahler and
let $H^0(X,\Omega^1_{X})$ denote the space of holomorphic one-forms
on $X$.  For every path $\gamma: I \ra  X$, 
$ \omega \mapsto \int_{\gamma} \omega$ defines a
linear map $H^0(X,\Omega^1_{X}) \ra \bbC$. This
correspondence induces a natural map 
$$\iota:  H_{1}(X,\bbZ) \ra H^0(X,\Omega^1_{X})^* \; . $$
Let us put $\Gamma = \iota(H_{1}(X,\bbZ))$. By Hodge-theory, 
$\Gamma \cong H_{1}(X,\bbZ)\big/ {\rm Torsion}$ is a lattice 
in $H^0(X,\Omega^1_{X})^*$. 
Define the Albanese torus $\Alb(X)$ of $X$ as 
$$ \Alb(X) = H^0(X,\Omega^1_{X})^*\big/\,\Gamma . $$
Now fix a base-point  
$x_{0} \in X$. Integration along paths defines a natural holomorphic map,
the Albanese morphism $ \alb_{X}: X \ra  \Alb(X)$, via
$$  \alb_{X}: \, \;   x \, \mapsto \int_{x_{0}}^x \; \; \;  . $$
A salient feature of the Albanese construction is that it reduces
integration of holomorphic one-forms on $X$ to integration of
linear one-forms on $\Alb(X)$.  In fact, let $\omega \in H^0(X,\Omega^1_{X})$
be a holomorphic one-form.
By duality of vector spaces, $\omega$ defines a (linear) holomorphic
one-form $\bar{\omega} \in H^0(\Alb(X), \Omega^1_{\Alb(X)})$. Then the identity
$$  \int_{\gamma} \alb_{X}^{*} \bar{\omega} \; = \int_{\gamma} \omega$$
is easily verified.  
The following lemma is a standard fact, and it is a direct consequence of the
preceding remark:

\begin{lemma} \label{lemma:albH1} 
The Albanese morphism $ \alb_{X}: X \ra  \Alb(X)$ 
induces an isomorphism of cohomology groups $ \alb_{X}^*: H^1(\Alb(X),\bbZ) \ra 
H^1(X,\bbZ)$.
\end{lemma}
\begin{proof}
Consider $\Gamma^*  \leq H^0(X, \Omega^1_{X})$, the dual lattice for $\Gamma$. 
Then the corresponding linear forms $\bar{\omega}$, $\omega \in \Gamma^*$, 
define one-forms in $\Omega^1(\Alb(X))$ which represent the 
integral cohomology $H^1(\Alb(X), \bbZ)$. Recall that $H^1(X, \bbZ)$ 
has no torsion. Then, by the definition of $\Gamma$, 
integration of forms on the cycles of $X$ 
identifies $\Gamma^*$ with $H^1(X, \bbZ)$.  Thus, the identity 
$ \alb_{X}^{*} \bar{\omega} = \omega$ proves the lemma.
\end{proof}

\subsection{The mapping degree of $\alb_{X}$}
We employ here the differential topological notion of
the mapping degree.  Let $f: M \ra N$ be a 
smooth map, where $M$ and $N$ are compact oriented smooth manifolds,
$k= \dim M = \dim N$. Let $y \in N$ be a regular value for $f$, then 
$$\deg f =  \sum_{x \in f^{-1}(y)}\! \!{\rm sign} \det(df_{x}) $$
is an integer which is independent of the choice of $y$, see \cite{Milnor_1}.  
The degree of $f$ also coincides with the unique integer which
is defined by the induced map $f^*: H^k(N,\bbZ) \ra H^k(M,\bbZ)$
relative to the orientation classes  of $M$ and $N$. 
Note the following particular useful consequence of the
definition of degree: 
\emph{if $\deg f$ is non-zero then $f$ must be surjective}. 

We shall apply this concept to holomorphic maps $f: X \ra Y$ between
compact complex manifolds $X$ and $Y$.  Here $k= 2n$,
where $n= \dim X = \dim Y$.  Both
$X$ and $Y$ have natural orientations determined by their
complex atlas. Moreover, with respect to these orientations, ${\rm sign} \det(df_{x}) = 1$, for all 
regular points $x \in X$. Thus if $f$ is  holomorphic, its
degree $\deg f$ is a non-negative integer, and $\deg f$ 
counts precisely the number of elements in the 
fibers over the regular values of $f$. 

Next we provide a criterion for finiteness of the fibers
of a holomorphic map.  
A proper holomorphic map $f:X \ra Y$ is called finite if for all $y \in Y$, $f^{-1}(y)$
is a finite set. 
\begin{lemma} \label{lemma:fibers} 
Assume that $X$ is compact K\"ahler and that the K\"ahler class
$[ \omega ]$ of $X$ is contained in $f^* H^2(Y,\bbR)$ then $f$ is finite. 
\end{lemma}
\begin{proof} It is enough to prove that the connected components
of $f^{-1}(y)$ are points. Let $F$, $\dim F=p$ be such a component,
and consider $\omega_{p} = \omega_{F}^p$ the restriction of the
$p$-th power of the K\"ahler form $\omega$ on $F$. Then $\omega_{p}$ is 
a positive definite volume-form on the analytic subvariety $F$ of $X$. 
The assumption on the K\"ahler class implies that 
$\omega_{p}$ is cohomologous to a pull-back 
form $f^*(\beta)$, where $\beta \in \Omega^{2p}(Y)$. The integral
of forms on analytic varieties is well defined, and by a version 
Stoke's theorem descends to cohomology classes, see \cite{Lelong}
and \cite[Theorem 11.21]{Voisin}.  
Since the pull-back $f^*(\beta)$ has zero integral on $F$, this
also holds for $\omega_{p}$. Hence, $F$ must be a point.
\end{proof}

We  are ready for the proof of our main result on the Albanese
mapping.  

\begin{proof}[Proof of Theorem 5]
Our strategy is to prove that the Albanese mapping $\alb_{X}: X \ra \Alb(X)$ 
is a holomorphic branched covering (a so-called {\it analytic covering}, see 
\cite[Ch.\ 7, \S 2]{Grauert-Remmert} for the definition) 
using the assumptions on the cohomology. 
Since $X$ and $\Alb(X)$ are compact complex manifolds, this amounts 
to showing that  $\alb_{X}$ is surjective
and has finite fibers, see 
\cite[Ch.\ 9 \S 3.3]{Grauert-Remmert}.  

Note that, by Lemma \ref{lemma:albH1},  the Albanese mapping 
$\alb_{X}$ induces an isomorphism $H^1(\Alb(X),\bbZ) \, \ra H^1(X,\bbZ)$. 
This and the assumption that $H^1(X,\bbZ)$ generates the top cohomology 
$H^{2n}(X,\bbZ)$ 
of $X$ implies that $\alb_{X}^*:  H^{2n}(\Alb(X),\bbZ) \ra H^{2n}(X,\bbZ)$
is an isomorphism on the integral top-cohomology groups. 
Thus $\deg \alb_{X} = 1$. In particular, $\alb_{X}$ 
is also a surjective map.

As a general fact, if $X$ is compact K\"ahler and $f: X \ra Y$ is a 
surjective holomorphic map then  (using techniques similar to the proof of
Lemma \ref{lemma:fibers}) the induced map $f^*: H^*(Y,\bbR) \ra H^*(X,\bbR)$
is injective, see \cite[Lemme 2.4]{Campana} or \cite[Lemma 7.28]{Voisin}. 
In particular, $\alb_{X}^* : H^2(\Alb(X),\bbC) \ra H^2(X,\bbC)$ is injective
and, hence, an isomorphism, since we have assumed 
that $\dim H^2(X,\bbC) \leq \dim H^2(\Alb(X),\bbC)$. 
Applying Lemma \ref{lemma:fibers} we can deduce that 
$\alb_{X}$ is a finite map and, hence, an analytic covering. 

Since $\alb_{X}$ has degree 1, it is injective at regular points. 
But at the singular points the number of elements in the
fibers cannot increase. For a proof of the latter fact, see for example 
\cite[Chapter 8, \S 1]{Grauert-Remmert}. Thus $\alb_{X}$ is bijective,
hence a holomorphic homeomorphism. 

The latter fact implies that $\alb_{X}$
is biholomorphic. This can be deduced as follows: 
The continuous map 
$\alb_{X}^{-1} : \Alb(X) \rightarrow X$ 
is holomorphic outside the set $S_{ \Alb(X)}\subset  \Alb(X)$ 
of singular values of  $\alb_{X}$. The set of singular points
$S_X$ of $\alb_{X}$ is clearly an analytic subset of $X$ 
and is mapped onto
$S_{ \Alb(X)}$ by the finite holomorphic map $alb_{X}$. This 
implies that $S_{ \Alb(X)}\subset  \Alb(X)$ is an analytic subset, by 
\cite[Proposition 4.1.6]{Grauert-Remmert}. 
Since $\Alb(X)$ is a complex manifold, the Riemann extension theorem 
\cite[Theorem 7.1.1]{Grauert-Remmert} shows that 
$\alb_{X}^{-1}$ is holomorphic. Thus  $\alb_{X}$ is biholomorphic.  
\end{proof}

% \section{Remaining proofs}

\end{document}